%%%%%%%%%%%%%%%%%%%%%%%%%%%%%%% including smallfonts.tex %%%%%%%%%%%%%%%%%%%%%%%%%%%%%%%
%smallfonts.tex
%
\newskip\ttglue
\font\fiverm=cmr5
\font\fivei=cmmi5
\font\fivesy=cmsy5
\font\fivebf=cmbx5
\font\sixrm=cmr6
\font\sixi=cmmi6
\font\sixsy=cmsy6
\font\sixbf=cmbx6
\font\sevenrm=cmr7
\font\eightrm=cmr8
\font\eighti=cmmi8
\font\eightsy=cmsy8
\font\eightit=cmti8
\font\eightsl=cmsl8
\font\eighttt=cmtt8
\font\eightbf=cmbx8
\font\ninerm=cmr9
\font\ninei=cmmi9
\font\ninesy=cmsy9
\font\nineit=cmti9
\font\ninesl=cmsl9
\font\ninett=cmtt9
\font\ninebf=cmbx9
\font\twelverm=cmr12
\font\twelvei=cmmi12
\font\twelvesy=cmsy12
\font\twelveit=cmti12
\font\twelvesl=cmsl12
\font\twelvett=cmtt12
\font\twelvebf=cmbx12

%% EIGHT POINT FONT FAMILY

\def\eightpoint{\def\rm{\fam0\eightrm}  
  \textfont0=\eightrm \scriptfont0=\sixrm \scriptscriptfont0=\fiverm
  \textfont1=\eighti  \scriptfont1=\sixi  \scriptscriptfont1=\fivei
  \textfont2=\eightsy  \scriptfont2=\sixsy  \scriptscriptfont2=\fivesy
  \textfont3=\tenex  \scriptfont3=\tenex  \scriptscriptfont3=\tenex
  \textfont\itfam=\eightit  \def\it{\fam\itfam\eightit}
  \textfont\slfam=\eightsl  \def\sl{\fam\slfam\eightsl}
  \textfont\ttfam=\eighttt  \def\tt{\fam\ttfam\eighttt}
  \textfont\bffam=\eightbf  \scriptfont\bffam=\sixbf
    \scriptscriptfont\bffam=\fivebf  \def\bf{\fam\bffam\eightbf}
  \tt  \ttglue=.5em plus.25em minus.15em
  \normalbaselineskip=9pt
  \setbox\strutbox=\hbox{\vrule height7pt depth2pt width0pt}
  \let\sc=\sixrm  \let\big=\eightbig \normalbaselines\rm}

\def\eightbig#1{{\hbox{$\textfont0=\ninerm\textfont2=\ninesy
        \left#1\vbox to6.5pt{}\right.$}}}

%% NINE POINT FONT FAMILY

\def\ninepoint{\def\rm{\fam0\ninerm}  
  \textfont0=\ninerm \scriptfont0=\sixrm \scriptscriptfont0=\fiverm
  \textfont1=\ninei  \scriptfont1=\sixi  \scriptscriptfont1=\fivei
  \textfont2=\ninesy  \scriptfont2=\sixsy  \scriptscriptfont2=\fivesy
  \textfont3=\tenex  \scriptfont3=\tenex  \scriptscriptfont3=\tenex
  \textfont\itfam=\nineit  \def\it{\fam\itfam\nineit}
  \textfont\slfam=\ninesl  \def\sl{\fam\slfam\ninesl}
  \textfont\ttfam=\ninett  \def\tt{\fam\ttfam\ninett}
  \textfont\bffam=\ninebf  \scriptfont\bffam=\sixbf
    \scriptscriptfont\bffam=\fivebf  \def\bf{\fam\bffam\ninebf}
  \tt  \ttglue=.5em plus.25em minus.15em
  \normalbaselineskip=11pt
  \setbox\strutbox=\hbox{\vrule height8pt depth3pt width0pt}
  \let\sc=\sevenrm  \let\big=\ninebig \normalbaselines\rm}

\def\ninebig#1{{\hbox{$\textfont0=\tenrm\textfont2=\tensy
        \left#1\vbox to7.25pt{}\right.$}}}

%% TWELVE POINT FONT FAMILY --- not really small

\def\twelvepoint{\def\rm{\fam0\twelverm}  
  \textfont0=\twelverm \scriptfont0=\eightrm \scriptscriptfont0=\sixrm
  \textfont1=\twelvei  \scriptfont1=\eighti  \scriptscriptfont1=\sixi
  \textfont2=\twelvesy  \scriptfont2=\eightsy  \scriptscriptfont2=\sixsy
  \textfont3=\tenex  \scriptfont3=\tenex  \scriptscriptfont3=\tenex
  \textfont\itfam=\twelveit  \def\it{\fam\itfam\twelveit}
  \textfont\slfam=\twelvesl  \def\sl{\fam\slfam\twelvesl}
  \textfont\ttfam=\twelvett  \def\tt{\fam\ttfam\twelvett}
  \textfont\bffam=\twelvebf  \scriptfont\bffam=\eightbf
    \scriptscriptfont\bffam=\sixbf  \def\bf{\fam\bffam\twelvebf}
  \tt  \ttglue=.5em plus.25em minus.15em
  \normalbaselineskip=11pt
  \setbox\strutbox=\hbox{\vrule height8pt depth3pt width0pt}
  \let\sc=\sevenrm  \let\big=\twelvebig \normalbaselines\rm}

\def\twelvebig#1{{\hbox{$\textfont0=\tenrm\textfont2=\tensy
        \left#1\vbox to7.25pt{}\right.$}}}
%%%%%%%%%%%%%%%%%%%%%%%%%%%%%%%%% end of smallfonts.tex %%%%%%%%%%%%%%%%%%%%%%%%%%%%%%%%
%%%%%%%%%%%%%%%%%%%%%%%%%%%%%%% including param.2 %%%%%%%%%%%%%%%%%%%%%%%%%%%%%%%
%param.2
\magnification=\magstep1
\def\firstpage{1}
\pageno=\firstpage
%%%%%%%%%%%%%%%%%%%%%%%%%%%%%%%%% end of param.2 %%%%%%%%%%%%%%%%%%%%%%%%%%%%%%%%
%%%%%%%%%%%%%%%%%%%%%%%%%%%%%%% including fonts.6 %%%%%%%%%%%%%%%%%%%%%%%%%%%%%%%
%fonts.6
\font\fiverm=cmr5
\font\sevenrm=cmr7
\font\sevenbf=cmbx7
\font\eightrm=cmr8
\font\eightbf=cmbx8
\font\ninerm=cmr9
\font\ninebf=cmbx9
\font\tenbf=cmbx10
\font\magtenbf=cmbx10 scaled\magstep1

\catcode`\@=11
%
%  Include all definitions related to the fonts msam, msbm and eufm, so that
%  when this file is used by itself, the results with respect to those fonts
%  are equivalent to what they would have been using AMS-TeX.
%  Most symbols in fonts msam and msbm are defined using \newsymbol;
%  however, a few symbols that replace composites defined in plain must be
%  defined with \mathchardef.

\def\undefine#1{\let#1\undefined}
\def\newsymbol#1#2#3#4#5{\let\next@\relax
 \ifnum#2=\@ne\let\next@\msafam@\else
 \ifnum#2=\tw@\let\next@\msbfam@\fi\fi
 \mathchardef#1="#3\next@#4#5}
\def\mathhexbox@#1#2#3{\relax
 \ifmmode\mathpalette{}{\m@th\mathchar"#1#2#3}%
 \else\leavevmode\hbox{$\m@th\mathchar"#1#2#3$}\fi}
\def\hexnumber@#1{\ifcase#1 0\or 1\or 2\or 3\or 4\or 5\or 6\or 7\or 8\or
 9\or A\or B\or C\or D\or E\or F\fi}

\font\tenmsa=msam10
\font\sevenmsa=msam7
\font\fivemsa=msam5
\newfam\msafam
\textfont\msafam=\tenmsa
\scriptfont\msafam=\sevenmsa
\scriptscriptfont\msafam=\fivemsa
\edef\msafam@{\hexnumber@\msafam}
\mathchardef\dabar@"0\msafam@39
\def\dashrightarrow{\mathrel{\dabar@\dabar@\mathchar"0\msafam@4B}}
\def\dashleftarrow{\mathrel{\mathchar"0\msafam@4C\dabar@\dabar@}}

\def\ulcorner{\delimiter"4\msafam@70\msafam@70 }
\def\urcorner{\delimiter"5\msafam@71\msafam@71 }
\def\llcorner{\delimiter"4\msafam@78\msafam@78 }
\def\lrcorner{\delimiter"5\msafam@79\msafam@79 }
%    Note that there should not be a final space after the digits for a
%    \mathhexbox@.
\def\yen{{\mathhexbox@\msafam@55}}
\def\checkmark{{\mathhexbox@\msafam@58}}
\def\circledR{{\mathhexbox@\msafam@72}}
\def\maltese{{\mathhexbox@\msafam@7A}}

\font\tenmsb=msbm10
\font\sevenmsb=msbm7
\font\fivemsb=msbm5
\newfam\msbfam
\textfont\msbfam=\tenmsb
\scriptfont\msbfam=\sevenmsb
\scriptscriptfont\msbfam=\fivemsb
\edef\msbfam@{\hexnumber@\msbfam}
\def\Bbb#1{{\fam\msbfam\relax#1}}
\def\widehat#1{\setbox\z@\hbox{$\m@th#1$}%
 \ifdim\wd\z@>\tw@ em\mathaccent"0\msbfam@5B{#1}%
 \else\mathaccent"0362{#1}\fi}
\def\widetilde#1{\setbox\z@\hbox{$\m@th#1$}%
 \ifdim\wd\z@>\tw@ em\mathaccent"0\msbfam@5D{#1}%
 \else\mathaccent"0365{#1}\fi}
\font\teneufm=eufm10
\font\seveneufm=eufm7
\font\fiveeufm=eufm5
\newfam\eufmfam
\textfont\eufmfam=\teneufm
\scriptfont\eufmfam=\seveneufm
\scriptscriptfont\eufmfam=\fiveeufm

\catcode`\@=11
%%  Load amssym.def if necessary: If \newsymbol is undefined, do nothing
%%  and the following \input statement will be executed; otherwise
%%  change \input to a temporary no-op.
%#\ifx\undefined\newsymbol \else \begingroup\def\input#1 {\endgroup}\fi
%#\input amssym.def \relax
%%  Most symbols in fonts msam and msbm are defined using \newsymbol.  A few
%%  that are delimiters or otherwise require special treatment have already
%%  been defined as soon as the fonts were loaded.  Finally, a few symbols
%%  that replace composites defined in plain must be undefined first.
\newsymbol\boxdot 1200
\newsymbol\boxplus 1201
\newsymbol\boxtimes 1202
\newsymbol\square 1003
\newsymbol\blacksquare 1004
\newsymbol\centerdot 1205
\newsymbol\lozenge 1006
\newsymbol\blacklozenge 1007
\newsymbol\circlearrowright 1308
\newsymbol\circlearrowleft 1309
\undefine\rightleftharpoons
\newsymbol\rightleftharpoons 130A
\newsymbol\leftrightharpoons 130B
\newsymbol\boxminus 120C
\newsymbol\Vdash 130D
\newsymbol\Vvdash 130E
\newsymbol\vDash 130F
\newsymbol\twoheadrightarrow 1310
\newsymbol\twoheadleftarrow 1311
\newsymbol\leftleftarrows 1312
\newsymbol\rightrightarrows 1313
\newsymbol\upuparrows 1314
\newsymbol\downdownarrows 1315
\newsymbol\upharpoonright 1316
 
\newsymbol\downharpoonright 1317
\newsymbol\upharpoonleft 1318
\newsymbol\downharpoonleft 1319
\newsymbol\rightarrowtail 131A
\newsymbol\leftarrowtail 131B
\newsymbol\leftrightarrows 131C
\newsymbol\rightleftarrows 131D
\newsymbol\Lsh 131E
\newsymbol\Rsh 131F
\newsymbol\rightsquigarrow 1320
\newsymbol\leftrightsquigarrow 1321
\newsymbol\looparrowleft 1322
\newsymbol\looparrowright 1323
\newsymbol\circeq 1324
\newsymbol\succsim 1325
\newsymbol\gtrsim 1326
\newsymbol\gtrapprox 1327
\newsymbol\multimap 1328
\newsymbol\therefore 1329
\newsymbol\because 132A
\newsymbol\doteqdot 132B
 
\newsymbol\triangleq 132C
\newsymbol\precsim 132D
\newsymbol\lesssim 132E
\newsymbol\lessapprox 132F
\newsymbol\eqslantless 1330
\newsymbol\eqslantgtr 1331
\newsymbol\curlyeqprec 1332
\newsymbol\curlyeqsucc 1333
\newsymbol\preccurlyeq 1334
\newsymbol\leqq 1335
\newsymbol\leqslant 1336
\newsymbol\lessgtr 1337
\newsymbol\backprime 1038
\newsymbol\risingdotseq 133A
\newsymbol\fallingdotseq 133B
\newsymbol\succcurlyeq 133C
\newsymbol\geqq 133D
\newsymbol\geqslant 133E
\newsymbol\gtrless 133F
\newsymbol\sqsubset 1340
\newsymbol\sqsupset 1341
\newsymbol\vartriangleright 1342
\newsymbol\vartriangleleft 1343
\newsymbol\trianglerighteq 1344
\newsymbol\trianglelefteq 1345
\newsymbol\bigstar 1046
\newsymbol\between 1347
\newsymbol\blacktriangledown 1048
\newsymbol\blacktriangleright 1349
\newsymbol\blacktriangleleft 134A
\newsymbol\vartriangle 134D
\newsymbol\blacktriangle 104E
\newsymbol\triangledown 104F
\newsymbol\eqcirc 1350
\newsymbol\lesseqgtr 1351
\newsymbol\gtreqless 1352
\newsymbol\lesseqqgtr 1353
\newsymbol\gtreqqless 1354
\newsymbol\Rrightarrow 1356
\newsymbol\Lleftarrow 1357
\newsymbol\veebar 1259
\newsymbol\barwedge 125A
\newsymbol\doublebarwedge 125B
\undefine\angle
\newsymbol\angle 105C
\newsymbol\measuredangle 105D
\newsymbol\sphericalangle 105E
\newsymbol\varpropto 135F
\newsymbol\smallsmile 1360
\newsymbol\smallfrown 1361
\newsymbol\Subset 1362
\newsymbol\Supset 1363
\newsymbol\Cup 1264
 
\newsymbol\Cap 1265
 
\newsymbol\curlywedge 1266
\newsymbol\curlyvee 1267
\newsymbol\leftthreetimes 1268
\newsymbol\rightthreetimes 1269
\newsymbol\subseteqq 136A
\newsymbol\supseteqq 136B
\newsymbol\bumpeq 136C
\newsymbol\Bumpeq 136D
\newsymbol\lll 136E
 
\newsymbol\ggg 136F
 
\newsymbol\circledS 1073
\newsymbol\pitchfork 1374
\newsymbol\dotplus 1275
\newsymbol\backsim 1376
\newsymbol\backsimeq 1377
\newsymbol\complement 107B
\newsymbol\intercal 127C
\newsymbol\circledcirc 127D
\newsymbol\circledast 127E
\newsymbol\circleddash 127F
\newsymbol\lvertneqq 2300
\newsymbol\gvertneqq 2301
\newsymbol\nleq 2302
\newsymbol\ngeq 2303
\newsymbol\nless 2304
\newsymbol\ngtr 2305
\newsymbol\nprec 2306
\newsymbol\nsucc 2307
\newsymbol\lneqq 2308
\newsymbol\gneqq 2309
\newsymbol\nleqslant 230A
\newsymbol\ngeqslant 230B
\newsymbol\lneq 230C
\newsymbol\gneq 230D
\newsymbol\npreceq 230E
\newsymbol\nsucceq 230F
\newsymbol\precnsim 2310
\newsymbol\succnsim 2311
\newsymbol\lnsim 2312
\newsymbol\gnsim 2313
\newsymbol\nleqq 2314
\newsymbol\ngeqq 2315
\newsymbol\precneqq 2316
\newsymbol\succneqq 2317
\newsymbol\precnapprox 2318
\newsymbol\succnapprox 2319
\newsymbol\lnapprox 231A
\newsymbol\gnapprox 231B
\newsymbol\nsim 231C
\newsymbol\ncong 231D
\newsymbol\diagup 201E
\newsymbol\diagdown 201F
\newsymbol\varsubsetneq 2320
\newsymbol\varsupsetneq 2321
\newsymbol\nsubseteqq 2322
\newsymbol\nsupseteqq 2323
\newsymbol\subsetneqq 2324
\newsymbol\supsetneqq 2325
\newsymbol\varsubsetneqq 2326
\newsymbol\varsupsetneqq 2327
\newsymbol\subsetneq 2328
\newsymbol\supsetneq 2329
\newsymbol\nsubseteq 232A
\newsymbol\nsupseteq 232B
\newsymbol\nparallel 232C
\newsymbol\nmid 232D
\newsymbol\nshortmid 232E
\newsymbol\nshortparallel 232F
\newsymbol\nvdash 2330
\newsymbol\nVdash 2331
\newsymbol\nvDash 2332
\newsymbol\nVDash 2333
\newsymbol\ntrianglerighteq 2334
\newsymbol\ntrianglelefteq 2335
\newsymbol\ntriangleleft 2336
\newsymbol\ntriangleright 2337
\newsymbol\nleftarrow 2338
\newsymbol\nrightarrow 2339
\newsymbol\nLeftarrow 233A
\newsymbol\nRightarrow 233B
\newsymbol\nLeftrightarrow 233C
\newsymbol\nleftrightarrow 233D
\newsymbol\divideontimes 223E
\newsymbol\varnothing 203F
\newsymbol\nexists 2040
\newsymbol\Finv 2060
\newsymbol\Game 2061
\newsymbol\mho 2066
\newsymbol\eth 2067
\newsymbol\eqsim 2368
\newsymbol\beth 2069
\newsymbol\gimel 206A
\newsymbol\daleth 206B
\newsymbol\lessdot 236C
\newsymbol\gtrdot 236D
\newsymbol\ltimes 226E
\newsymbol\rtimes 226F
\newsymbol\shortmid 2370
\newsymbol\shortparallel 2371
\newsymbol\smallsetminus 2272
\newsymbol\thicksim 2373
\newsymbol\thickapprox 2374
\newsymbol\approxeq 2375
\newsymbol\succapprox 2376
\newsymbol\precapprox 2377
\newsymbol\curvearrowleft 2378
\newsymbol\curvearrowright 2379
\newsymbol\digamma 207A
\newsymbol\varkappa 207B
\newsymbol\Bbbk 207C
\newsymbol\hslash 207D
\undefine\hbar
\newsymbol\hbar 207E
\newsymbol\backepsilon 237F
%  Restore the catcode value for @ that was previously saved.
%#\catcode`\@=\csname pre amssym.tex at\endcsname

%\endinput
%%%%%%%%%%%%%%%%%%%%%%%%%%%%%%%%% end of symbols.1 %%%%%%%%%%%%%%%%%%%%%%%%%%%%%%%%
%%%%%%%%%%%%%%%%%%%%%%%%%%%%%%% including links.1 %%%%%%%%%%%%%%%%%%%%%%%%%%%%%%%
% links.1
% adapted from http://insti.physics.sunysb.edu/~siegel/tex.shtml
%
% postscript/pdf
\newcount\marknumber	\marknumber=1
\newcount\countdp \newcount\countwd \newcount\countht 
%
% for ordinary tex
%
\ifx\pdfoutput\undefined
\def\rgboo#1{}
\def\postscript#1{\special{" #1}}		% for dvips
\postscript{
	/bd {bind def} bind def
	/fsd {findfont exch scalefont def} bd
	/sms {setfont moveto show} bd
	/ms {moveto show} bd
	/pdfmark where		% printers ignore pdfmarks
	{pop} {userdict /pdfmark /cleartomark load put} ifelse
	[ /PageMode /UseOutlines		% bookmark window open
	/DOCVIEW pdfmark}
\def\bookmark#1#2{\postscript{		% #1=subheadings (if not 0)
	[ /Dest /MyDest\the\marknumber /View [ /XYZ null null null ] /DEST pdfmark
	[ /Title (#2) /Count #1 /Dest /MyDest\the\marknumber /OUT pdfmark}%
	\advance\marknumber by1}
\def\pdfclink#1#2#3{%
	\hskip-.25em\setbox0=\hbox{#2}%
		\countdp=\dp0 \countwd=\wd0 \countht=\ht0%
		\divide\countdp by65536 \divide\countwd by65536%
			\divide\countht by65536%
		\advance\countdp by1 \advance\countwd by1%
			\advance\countht by1%
		\def\linkdp{\the\countdp} \def\linkwd{\the\countwd}%
			\def\linkht{\the\countht}%
	\postscript{
		[ /Rect [ -1.5 -\linkdp.0 0\linkwd.0 0\linkht.5 ] 
		/Border [ 0 0 0 ]
		/Action << /Subtype /URI /URI (#3) >>
		/Subtype /Link
		/ANN pdfmark}{\rgb{#1}{#2}}}
%
% for pdftex
%
\else
\def\rgboo#1{\pdfliteral{#1 rg #1 RG}}
\pdfcatalog{/PageMode /UseOutlines}		% bookmark window open
\def\bookmark#1#2{
	\pdfdest num \marknumber xyz
	\pdfoutline goto num \marknumber count #1 {#2}
	\advance\marknumber by1}
\def\pdfklink#1#2{%
	\noindent\pdfstartlink user
		{/Subtype /Link
		/Border [ 0 0 0 ]
		/A << /S /URI /URI (#2) >>}{\rgb{1 0 0}{#1}}%
	\pdfendlink}
\fi

\def\rgbo#1#2{\rgboo{#1}#2\rgboo{0 0 0}}
\def\rgb#1#2{\mark{#1}\rgbo{#1}{#2}\mark{0 0 0}}
\def\pdfklink#1#2{\pdfclink{1 0 0}{#1}{#2}}
\def\pdflink#1{\pdfklink{#1}{#1}}
%
% examples:
% \bookmark{0}{look here}
% \pdfclink{0 0 1}{testlink}{http://www.google.com/}
% \pdfklink{testlink}{http://www.google.com/}
% \pdflink{http://www.google.com/}
%%%%%%%%%%%%%%%%%%%%%%%%%%%%%%%%% end of links.1 %%%%%%%%%%%%%%%%%%%%%%%%%%%%%%%%
%%%%%%%%%%%%%%%%%%%%%%%%%%%%%%% including titles.8c %%%%%%%%%%%%%%%%%%%%%%%%%%%%%%%
%titles.8
% requires fonts.5 or higher and smallfonts.tex
% uses links.* if included
% enumerates \demo consecutively (no section number)
%
\newcount\seccount  %% sections
\newcount\subcount  %% subsection
\newcount\clmcount  %% claim
\newcount\equcount  %% equation
\newcount\refcount  %% reference
\newcount\demcount  %% example
\newcount\execount  %% exercise
\newcount\procount  %% problem
\seccount=0
\equcount=1
\clmcount=1
\subcount=1
\refcount=1
\demcount=0
\execount=0
\procount=0
%
%% MISC STUFF
\def\proof{\medskip\noindent{\bf Proof.\ }}
\def\proofof(#1){\medskip\noindent{\bf Proof of \csname c#1\endcsname.\ }}
\def\qed{\hfill{\sevenbf QED}\par\medskip}
\def\references{\bigskip\noindent\hbox{\bf References}\medskip
                \ifx\pdflink\undefined\else\bookmark{0}{References}\fi}
\def\addref#1{\expandafter\xdef\csname r#1\endcsname{\number\refcount}
    \global\advance\refcount by 1}

\def\nextremark #1\par{\item{$\circ$} #1}
\def\firstremark #1\par{\bigskip\noindent{\bf Remarks.}
     \smallskip\nextremark #1\par}
\def\abstract#1\par{{\baselineskip=10pt
    \eightpoint\narrower\noindent{\eightbf Abstract.} #1\par}}
%
%% EQUATION
\def\equtag#1{\expandafter\xdef\csname e#1\endcsname{(\number\seccount.\number\equcount)}
              \global\advance\equcount by 1}
\def\equation(#1){\equtag{#1}\eqno\csname e#1\endcsname}
\def\equ(#1){\hskip-0.03em\csname e#1\endcsname}
%
%% CLAIMS (theorems etc)
\def\clmtag#1#2{\expandafter\xdef\csname c#2\endcsname{#1~\number\seccount.\number\clmcount}
                \global\advance\clmcount by 1}
\def\claim #1(#2) #3\par{\clmtag{#1}{#2}
    \vskip.1in\medbreak\noindent
    {\bf \csname c#2\endcsname .\ }{\sl #3}\par
    \ifdim\lastskip<\medskipamount
    \removelastskip\penalty55\medskip\fi}
\def\clm(#1){\csname c#1\endcsname}
%
%% SECTION
\def\sectag#1{\global\advance\seccount by 1
              \expandafter\xdef\csname sectionname\endcsname{\number\seccount. #1}
              \equcount=1 \clmcount=1 \subcount=1 \execount=0 \procount=0}
\def\section#1\par{\vskip0pt plus.1\vsize\penalty-40
    \vskip0pt plus -.1\vsize\bigskip\bigskip
    \sectag{#1}
    \message{\sectionname}\leftline{\magtenbf\sectionname}
    \nobreak\smallskip\noindent
    \ifx\pdflink\undefined
    \else
      \bookmark{0}{\sectionname}
    \fi}
%
%% SUBSECTION
\def\subtag#1{\expandafter\xdef\csname subsectionname\endcsname{\number\seccount.\number\subcount. #1}
              \global\advance\subcount by 1}
\def\subsection#1\par{\vskip0pt plus.05\vsize\penalty-20
    \vskip0pt plus -.05\vsize\medskip\medskip
    \subtag{#1}
    \message{\subsectionname}\leftline{\tenbf\subsectionname}
    \nobreak\smallskip\noindent
    \ifx\pdflink\undefined
    \else
      \bookmark{0}{.... \subsectionname}  %% can get a bit cluttered
    \fi}
%
%% DEMO (examples etc)
\def\demtag#1#2{\global\advance\demcount by 1
              \expandafter\xdef\csname de#2\endcsname{#1~\number\demcount}}
\def\demo #1(#2) #3\par{
  \demtag{#1}{#2}
  \vskip.1in\medbreak\noindent
  {\bf #1 \number\demcount.\enspace}
  {\rm #3}\par
  \ifdim\lastskip<\medskipamount
  \removelastskip\penalty55\medskip\fi}
\def\dem(#1){\csname de#1\endcsname}
%
%% EXERCISE
\def\exetag#1{\global\advance\execount by 1
              \expandafter\xdef\csname ex#1\endcsname{Exercise~\number\seccount.\number\execount}}
\def\exercise(#1) #2\par{
  \exetag{#1}
  \vskip.1in\medbreak\noindent
  {\bf Exercise \number\execount.}
  {\rm #2}\par
  \ifdim\lastskip<\medskipamount
  \removelastskip\penalty55\medskip\fi}
\def\exe(#1){\csname ex#1\endcsname}
%
%% PROBLEM
\def\protag#1{\global\advance\procount by 1
              \expandafter\xdef\csname pr#1\endcsname{\number\seccount.\number\procount}}
\def\problem(#1) #2\par{
  \ifnum\procount=0
    \parskip=6pt
    \vbox{\bigskip\centerline{\bf Problems \number\seccount}\nobreak\medskip}
  \fi
  \protag{#1}
  \item{\number\procount.} #2}
\def\pro(#1){Problem \csname pr#1\endcsname}
%
%%%%%%%%%%%%%%%%%%%%%%%%%%%%%%%%% end of titles.8c %%%%%%%%%%%%%%%%%%%%%%%%%%%%%%%%
%%%%%%%%%%%%%%%%%%%%%%%%%%%%%%% including macros.19 %%%%%%%%%%%%%%%%%%%%%%%%%%%%%%%
%macros.19
% requires fonts.5 or later
\def\rightheadline{\hfil}
\def\leftheadline{\sevenrm\hfil HANS KOCH\hfil}
\headline={\ifnum\pageno=\firstpage\hfil\else
\ifodd\pageno{{\fiverm\rightheadline}\number\pageno}
\else{\number\pageno\fiverm\leftheadline}\fi\fi}
\footline={\ifnum\pageno=\firstpage\hss\tenrm\folio\hss\else\hss\fi}

\let\cl=\centerline

\let\sss=\scriptscriptstyle

\def\natural{{\Bbb N}}

\def\real{{\Bbb R}}
\def\complex{{\Bbb C}}

\def\iso{{\Bbb J}}

\def\id{{\rm I}}

\def\bdot{\hbox{\bf .}}

\def\half{{1\over 2}}

\def\thalf{{\textstyle\half}}

\def\twomat#1#2#3#4{\left[\matrix{#1&#2\cr#3&#4\cr}\right]}

%

%

%
% from TeX book: used for commutative diagram
% in math mode, before using matrix, do
% \def\normalbaselines{\baselineskip20pt\lineskip3pt\lineskiplimit3pt}

%
\def\AA{{\cal A}}
\def\BB{{\cal B}}

\def\EE{{\cal E}}
\def\FF{{\cal F}}
\def\GG{{\cal G}}

\def\KK{{\cal K}}

\def\OO{{\cal O}}

\def\RR{{\cal R}}

\def\VV{{\cal V}}

\def\ZZ{{\cal Z}}
%

%%%%%%%%%%%%%%%%%%%%%%%%%%%%%%%%% end of macros.19 %%%%%%%%%%%%%%%%%%%%%%%%%%%%%%%%
%%%%%%%%%%%%%%%%%%%%%%%%%%%%% including sketch11.tex %%%%%%%%%%%%%%%%%%%%%%%%%%%%%
%\input smallfonts.tex
%\input param.2
%\input fonts.6
%\input symbols.1
%\input links.1
%\input titles.8c
%\input macros.19
%%%%%%%%%%%%%%%%%%%%%%%%%%%%%%%%%%%%%% some more symbols
\let\ss=\scriptstyle
\def\ssreal{{\ss\real}}
\def\transpose{{\sss\top}}
\def\ssperp{{\sss\perp}}
\def\ssold{{\sss\rm old}}

\def\rmO{{\rm O}}
\def\rmU{{\rm U}}
\def\iso{{\rm J}}
\def\range{{\rm range}}

\def\GL{{\rm GL}}
\def\Sp{{\rm Sp}}
\def\sp{{\rm sp}}
\def\Hess{{\Bbb H}}
\def\Kess{{\Bbb K}}
\def\llb{{[\![}}
\def\rrb{{]\!]}}
%%%%%%%%%%%%%%%%%%%%%%%%%%%%%%%%%%%%%% references
\addref{Jung}
\addref{DragtFinn}
\addref{Helgason}
\addref{DragtAbell}
\addref{Moser}
\addref{FengWang}
\addref{BE}
\addref{Lomeli}
\addref{Forest}
\addref{Bondt}
\addref{McLachlan}
\addref{KL}
%%%%%%%%%%%%%%%%%%%%%%%%%%%%%%%%%%%%%%%%%%%%%%%%%%%%%%%%%%%%%%%%%
\def\rightheadline{\sevenrm\hfil H.~Koch and H.E.~Lomel\'\i\hfil}
\def\leftheadline{\sevenrm\hfil On Hamiltonian flows
whose orbits are straight lines\hfil}
\cl{{\magtenbf On Hamiltonian flows whose orbits are straight lines}}
\bigskip

\cl{Hans Koch
\footnote{$^1$}
{{\sevenrm Department of Mathematics, University of Texas at Austin,
Austin, TX 78712}}
and H\'ector E.~Lomel\'\i\ 
\footnote{$^2$}
{{\sevenrm Department of Mathematics, Instituto Tecnol\'ogico Aut\'onomo
de M\'exico, M\'exico DF 01000, Mexico}}
\footnote{}
{\vskip-15pt{\sevenrm Current address: Department of Mathematics,
University of Texas at Austin, Austin, TX 78712}}}

\bigskip
\abstract
We consider real analytic Hamiltonians on $\ssreal^n\times\ssreal^n$
whose flow depends linearly on time.
Trivial examples are Hamiltonians $H(q,p)$ that
do not depend on the coordinate $q\in\ssreal^n$.
By a theorem of Moser [\rMoser], every polynomial Hamiltonian
of degree $3$ reduces to such a $q$-independent Hamiltonian
via a linear symplectic change of variables.
We show that such a reduction is impossible, in general,
for polynomials of degree $4$ or higher.
But we give a condition that implies linear-symplectic conjugacy
to another simple class of Hamiltonians.
The condition is shown to hold for all nondegenerate Hamiltonians
that are homogeneous of degree $4$.

\section Introduction and main results
%%%%%%%%%%%%%%%%%%%%%%%%%%%%%%%%%%%%%%

Polynomial Hamiltonians and maps have been studied extensively
and for a variety of different reasons.
Among other things, they constitute local normal forms
for more general Hamiltonians and maps,
and they provide a convenient testing ground
for new ideas in dynamical systems.
The restriction to polynomials also adds interesting algebraic
aspects to the problem.
This includes the possibility of classifying
polynomial maps with a given property,
and of decomposing them into simpler ones.

The work presented here was motivated by a question,
described below, that was raised in [\rLomeli] concerning symplectic
maps $F$ with $F-\id$ a homogeneous polynomial.
A differentiable map $F:\real^{2n}\to\real^{2n}$ is said to be symplectic if
$$
DF(x)^\transpose\iso DF(x)=\iso\,,
\qquad \iso=\twomat{0}{\id}{-\id}{0}\,,
\qquad x\in\real^{2n}\,.
\equation(FSymplectic)
$$
Here $DF(x)$ denotes the derivative of $F$ at $x$,
and $DF(x)^\transpose$ denotes its transpose (as a matrix).
If $F-\id$ is a homogeneous polynomial,
then it can be shown [\rKL] that $F$
is the time-one map of a Hamiltonian $H$,
and that $H$ is affine-integrable as defined below.
Thus, it is natural in this context to work with Hamiltonian systems.
Moreover, some of our results do not require
that the Hamiltonian be homogeneous, or even polynomial.

Let $H$ be a smooth function on $\real^{2n}$.
One of the basic facts from Hamiltonian mechanics
is that the vector field $X=\iso\nabla H$
defines a flow $\Phi:(t,x)\mapsto\Phi^t(x)$
whose time-$t$ maps $\Phi^t$ are symplectic.

\claim Definition(Jolt)
We say that a Hamiltonian $H$ is affine-integrable
if its flow $\Phi$ is linear in time:
$$
\Phi^t=\id+tX\,,\qquad X=\iso\nabla H\,,\qquad t\in\real\,.
\equation(Jolt)
$$

As we will see later, a Hamiltonian is affine-integrable
if and only if the corresponding vector field $X$
is constant along each orbit.
That is, for all $t\in\real$ we have
$$
X\circ(\id+tX)=X\,.
\equation(XoItX)
$$

Polynomial maps $F:\real^m\to\real^m$ with the property
that $X=F-\id$ satisfies \equ(XoItX)
are also called quasi-translations.
They arise naturally in the study of singular Hessians [\rBE].
To be more precise,
the standard definition of a quasi-translation
only requires that the identity \equ(XoItX) be satisfied for $t=1$.
However, this identity extends to $t\in\natural$ by induction,
and further to $t\in\complex$ if $X$ is a polynomial
(since $X\circ(\id+tX)-X$ is a polynomial in $t$
with infinitely many zeros).
So $\id+tX$ is the time-$t$ map for the vector field $X$.
Differentiating the identity \equ(XoItX) with respect to $t$ yields $(DX)X=0$.
Or equivalently, $(X^\transpose\nabla)^2\ell=0$
for all linear functions $\ell$.
This ``local nilpotency'' property is an alternative
way of characterizing quasi-translations [\rBondt]
and affine-integrable Hamiltonians [\rDragtAbell,\rFengWang].

In numerical analysis and physics,
symplectic quasi-translations are also called {\sl jolt maps}.
They constitute the basic building blocks in the so-called
Dragt-Finn factorization [\rDragtFinn] of more general symplectic maps.
This factorization has proved to be very useful
in symplectic numerical schemes, including the simulation of
Hamiltonian flows in plasmas [\rDragtAbell,\rForest].

{}From a dynamical systems point of view,
affine-integrable Hamiltonians are rather simple.
Not only is the vector field $X$ constant along each orbit,
but its components $X_j$ are Poisson-commuting invariants, as we will see later.
So an affine-integrable Hamiltonian $H$ is Liouville integrable,
at least if it satisfies a suitable nondegeneracy condition.
In addition, the geometry defined by the invariants $X_j$ is
quite restricted:

\claim Theorem(AffineSub)
Let $H$ be a real analytic affine-integrable Hamiltonian on $\real^{2n}$.
Then $H$ and its vector field $X$
are constant on the affine subspaces $x+\range(DX(x))$.
If $DX(x)$ has rank $n$ then $x+\range(DX(x))$
is a local level set for $X$.

The only affine-integrable Hamiltonians that we have been able to find
in the literature are all linear-symplectically conjugate
to Hamiltonians of the form $H(q,p)=K(p)$.
The time-$t$ map for such a Hamiltonian $H$
is a shear: $\Phi^t(q,p)=\bigl(q+t\nabla K(p),p\bigr)$.

\claim Definition(ShearHam)
We call $H:\real^{2n}\to\real$ a shear Hamiltonian
if $H$ is linear-symplectically conjugate to a Hamiltonian
that does not depend on the variable $q$.
In other words, there exists a linear symplectic change of variables
$U:\real^{2n}\to\real^{2n}$
such that $(H\circ U)(q,p)=K(p)$ for some function $K:\real^n\to\real$.

\demo Remark(Iwasawa)
As we will describe later, the linear map $U$ in
the above definition can be chosen both symplectic and orthogonal
(as a matrix).

\medskip
One of our goals is to find an affine-integrable Hamiltonian
that is not a shear, or to prove that there is no such Hamiltonian.
Partial non-existence results can be obtained by restricting
the class of Hamiltonians being considered.
A trivial case:
If $n=1$ then line-orbits are necessarily parallel,
so if $H:\real^2\to\real$ is affine-integrable, then there exists
a rotation $U$ of $\real^2$, such that $(H\circ U)(q,p)$ is independent of $q$.
In other cases the \clm(ShearHam) cannot be used directly.
We shall give now give an alternative characterization of shear Hamiltonians.
It applies to Hamiltonians that are regular
in the following sense:

\claim Definition(Regular)
We say that a Hamiltonian $H$ is regular if
there exists a point $x$ where $X(x)$ belongs to the
range of $DX(x)$.

Notice that homogeneous Hamiltonians of degree $\ge 2$
are regular, since their vector field vanishes at the origin.
The following theorem was proved in [\rLomeli]
for Hamiltonians that are homogeneous polynomials of degree $\ge 3$.

\claim Theorem(JoltEquiv)
A regular real analytic function $H$ on $\real^{2n}$
is a shear Hamiltonian
if and only if $DX(x)DX(y)=0$ holds for all $x,y\in\real^{2n}$.

To continue our discussion of special cases,
assume that $H$ is regular and affine-integrable.
Then $(DX)^2=0$, as we will see later.
In particular, if $H$ is quadratic then
\clm(JoltEquiv) implies that $H$ is a shear Hamiltonian.
The cubic case is covered by a result of Moser [\rMoser]
on quadratic symplectic maps on $\real^{2n}$.
It states that every such map $F$
admits a decomposition $F=A\circ S\circ L$
into three simple symplectic maps:
an affine map $A$, a shear $S(q,p)=(q+s(p),p)$, and a linear map $L$.
It is not hard to show that this result implies
-- and is essentially equivalent to -- the statement
that every homogeneous affine-integrable Hamiltonian
of degree $3$ is a shear.
A direct proof is given in Section 2.

This raises the question [\rLomeli] whether
every homogeneous affine-integrable Hamiltonian on $\real^{2n}$ is a shear.
Locally, much more is true:
It is well known that every smooth Hamiltonian $H$
is of the form $H(q,p)=K(p)$ in some local symplectic chart,
near any point where the vector field does not vanish.
The local conjugacy (chart) is nonlinear in general.
But if all orbits for $H$ are straight lines with constant velocity,
as is the case for affine-integrable Hamiltonians,
and if $X$ is constant on $n$-dimensional affine subspaces,
then one might think that this conjugacy can be chosen to be linear.
However, this is not true in general:

\claim Theorem(CounterExample)
The following Hamiltonian on $\real^8$ is affine-integrable but not a shear:
$$
H(q,p)
=q_1p_3^3+\sqrt{3}q_2p_3^2p_4+p_1p_4^3-\sqrt{3}p_2p_3p_4^2\,,
\qquad q,p\in\real^4\,.
\equation(CounterExample)
$$
In addition, $H$ is nondegenerate in the sense defined below.

\claim Definition(Nondegenerate)
A real analytic vector field $X$ on $\real^{2n}$
is said to be nondegenerate
if $DX(x)$ has rank $\ge n$ at some point $x\in\real^{2n}$.
If $X=\iso\nabla H$ then we also say that $H$ is nondegenerate.

We would like to stress that this is the one and only notion of
nondegeneracy that will be used in this paper.
Notice that, if $DX(x)^2=0$,
then the rank of $DX(x)$ can be no larger than $n$.
And if $X$ is analytic, then the rank is constant
outside some analytic set of codimension one.

The example \equ(CounterExample) belongs to a simple
class of Hamiltonians that we shall now describe.
Let $0\le d<n$. To simplify the description,
we write $q=(Q,\bar q)$ and $p=(P,\bar p)$,
where $Q,P\in\real^d$ and $\bar q,\bar p\in\real^{n-d}$.
Consider a Hamiltonian of the form
$$
H(q,p)=K(\bar p)+Q^\transpose V(\bar p)+P^\transpose W(\bar p)\,,
\equation(SimpleHam)
$$
with $K:\real^{n-d}\to\real$ and $V,W:\real^{n-d}\to\real^d$ differentiable.
Notice that $H$ does not depend on $\bar q$,
and thus $\bar p$ stays fixed under the flow.
Furthermore, the coordinates $Q$ and $P$
evolve linearly (in time) under the flow.
If $\bar q$ evolves linearly as well, then $H$ is affine-integrable.
As we will see later, this is the case if and only if
$$
W(\bar p)^\transpose DV(\bar p)-V(\bar p)^\transpose DW(\bar p)=0\,.
\equation(WDVVDWo)
$$
If $d=0$, then $P=Q=0$ and $\bar p$ can be identified with $p$.
In this case, \equ(SimpleHam) becomes $H(q,p)=K(p)$,
so $H$ is a shear Hamiltonian.

\demo Remark(HamDecomp)
The Hamiltonian \equ(SimpleHam) can be written as the sum of
$H_1=Q^\transpose V(\bar p)$ and $H_2=K(\bar p)+P^\transpose W(\bar p)$.
What makes this decomposition interesting is that
the Poisson bracket $\{H_1,H_2\}=(\nabla H_1)^\transpose\iso(\nabla H_2)$
of $H_1$ and $H_2$ Poisson-commutes with both $H_1$ and $H_2$.

\claim Theorem(SimpleHam)
Let $H$ be a nondegenerate real analytic affine-integrable Hamiltonian on $\real^{2n}$.
Then $H$ is linear-symplectically conjugate
to a Hamiltonian of the form \equ(SimpleHam)
if and only if $DX(x)DX(y)DX(z)=0$ for all $x,y,z\in\real^{2n}$.

Here, as in \clm(AffineSub),
a simple class of affine-integrable Hamiltonian is characterized
by a nilpotency-type condition on the derivative
of the vector field.
This suggest there may be a natural hierarchy of such conditions,
which characterize classes of increasingly complex affine-integrable Hamiltonians.

We will prove that the condition on $DX$ in \clm(SimpleHam)
holds if $H$ is a homogeneous polynomials of degree $4$.
As a result we obtain

\claim Theorem(SimplerHam)
Let $H$ be a nondegenerate affine-integrable Hamiltonian on $\real^{2n}$.
If $H$ is homogeneous of degree $4$
then $H$ is is linear-symplectically conjugate
to a Hamiltonian of the form \equ(SimpleHam).

As mentioned earlier, any symplectic map $F:\real^{2n}\to\real^{2n}$
with the property that $F-\id$ is a homogeneous polynomial
is the time-one map of an affine-integrable Hamiltonian [\rKL].
Combining this result with \clm(SimplerHam)
and using the decomposition described in \dem(HamDecomp),
we obtain the following factorization theorem.

\claim Theorem(CubicDecomp)
Let $F$ be a symplectic map on $\real^{2n}$
such that $F-\id$ is homogeneous of degree $3$ and nondegenerate.
Then $F$ admits a decomposition $F=F_1\circ F_4$,
where $F_1$ and $F_4$ are the time-one maps of two shear Hamiltonians.

\def\foo{\hskip 2pt}
\medskip
For the proofs of Theorems 1.2,\foo 1.5,\foo 1.6,\foo 1.8,\foo 1.9, and 1.10,
we refer to Sections 3,\foo 2,\foo 6,\foo 4,\foo 5, and 7, respectively.
Some additional results can be found in Section 2.

\section Basic properties
%%%%%%%%%%%%%%%%%%%%%%%%%

In the remaining part of this paper we always assume
that $H$ is a real analytic Hamiltonian on $\real^{2n}$.
Furthermore, by a ``homogeneous'' Hamiltonian we always mean
a homogeneous polynomial.

\medskip
As is true in general, the Hamiltonian $H$ is invariant under the flow
that it generates, so $(DH)X=0$.
Assuming that $H$ is affine-integrable, $X\circ(\id+tX)=X$, and thus $(DX)X=0$.
Furthermore, $\Phi^t=\id+tX$ is symplectic, which by \equ(FSymplectic) yields
$$
\iso+t\bigl[(DX)^\transpose\iso+\iso DX\bigr]
+t^2\bigl[(DX)^\transpose\iso DX\bigr]=\iso\,.
\equation(PhiIsoPhi)
$$
Using that the terms of order $t$ and $t^2$ have to vanish separately,
we get $(DX)^2=0$.
Differentiating the identity $X\circ(\id+tX)=X$
yields $[DX\circ(\id+tX)](\id+tDX)=DX$.
Multiplying on the right by $(\id-tDX)$ and using that $(DX)^2=0$,
we find in addition that $DX\circ(\id+tX)=DX$.
In summary, we have the following

\claim Lemma(HXDXinvariant)
Let $H$ be an affine-integrable Hamiltonian.
Then the functions $H$ and $X$ and $DX$ are constant along every orbit.
Furthermore, $(DH)X=0$ and $(DX)X=0$ and $(DX)^2=0$.

An equivalent formulation of these facts can be given in terms of Poisson brackets.
Assume that $H$ is affine-integrable.
Let $\ell(x)=(\iso u)^\transpose x$ for some vector $u\in\real^{2n}$.
Then $\{\ell,H\}$ is the directional derivative of $H$
in the direction $u$, which we denote by $\partial_uH$.
Being linear in the coordinate $x$, $\ell$ evolves linearly in time,
so $\{\ell,H\}=\partial_uH$ is invariant under the flow.
This implies the first identity in
$$
\{\partial_u H,H\}=0\,,\qquad
\{\partial_u\partial_v H,H\}=0\,,\qquad
\{\partial_u H,\partial_v H\}=0\,.
\equation(CommDuHDvH)
$$
The second and third identities are obtained from the first
by applying a derivative $\partial_v$. This yields
$\{\partial_u\partial_v H,H\}+\{\partial_u H,\partial_v H\}=0$,
and the two terms have to vanish separately
since the first is symmetric in $(u,v)$
and the second antisymmetric.
If $H$ is nondegenerate, then \equ(CommDuHDvH) shows that
$n$ of the vector field components $X_j$ constitute a
maximal set of Poisson-commuting invariants.
So $H$ is Liouville integrable, as mentioned earlier.

\smallskip
Let $x$ be a fixed but arbitrary point in $\real^{2n}$.
In the canonical splitting $\real^{2n}=\real^n\times\real^n$,
we can represent the derivative of $X$
and the Hessian $\Hess(x)=\nabla DH(x)$
as $2\times 2$ matrices
whose entries are $n\times n$ matrices,
$$
DX(x)=\twomat{\ZZ(x)^\transpose}{\AA(x)}{-\BB(x)}{-\ZZ(x)}\,,\qquad
\Hess(x)=\twomat{\BB(x)}{\ZZ(x)}{\ZZ(x)^\transpose}{\AA(x)}\,.
\equation(GeneralHess)
$$
Given that $X=\iso\nabla H$, we have $DX=\iso\Hess$.
Since $\Hess$ is symmetric, so are $\AA$ and $\BB$.
In the case of an affine-integrable Hamiltonian, $\AA\ZZ$ is symmetric as well,
as a result of the identity $(DX)^2=0$.

\claim Lemma(BZzero)
A regular Hamiltonian is of the form $H(q,p)=K(p)$
if and only if $\ZZ(x)=0$ and $\BB(x)=0$ for all $x$.

\proof
The necessity of the conditions $\ZZ=0$ and $\BB=0$ is obvious.
Assume now that they are satisfied.
Let $x_0=(q_0,p_0)$ be a point where
$X(x_0)=\bigl(\nabla_p H(x_0),-\nabla_qH(x_0)\bigr)$
belongs to the range of $DX(x_0)$. At this point we have $\nabla_q H(x_0)=0$.
Given that $D_q^2H=0$ by assumption, this implies
that the function $q\mapsto H(q,p_0)$ is constant.
Furthermore, $D_pH(q,p)$ does not depend on $q$, since $D_qD_pH=0$.
So $H(q,p)$ is independent of $q$ as well, since
$H(q,p)=H(q,p_0)+\int_0^1 D_pH(q,p_0+sv)v\,ds$ with $v=p-p_0$.
\qed

\claim Lemma(DXNormalForm)
Let $H$ be an affine-integrable Hamiltonian.
Given $x\in\real^{2n}$,
there exists an orthogonal symplectic $2n\times 2n$ matrix $U$,
and a diagonal $n\times n$ matrix $A$, such that
$$
U^{-1}DX(x)U=\twomat{0}{A}{0}{0}\,,\qquad
U^\transpose\Hess(x)U=\twomat{0}{0}{0}{A}\,.
\equation(DXNormalForm)
$$

\proof
Let $d$ be the rank of $M=\Hess(x)$.
Let $(u_1,u_2,\ldots,u_d)$ be an orthonormal set of
eigenvectors for the nonzero eigenvalues of $M$.
Since $M\iso M=0$,
the vectors $\iso u_j$ are eigenvectors of $M$
for the eigenvalue $0$.
Consider first the case $d=n$.
Let $U$ be the $2n\times 2n$ matrix whose columns vectors
are $\iso u_1,\ldots,\iso u_n,u_1,\ldots,u_n$, in this order.
Clearly, $U$ is orthogonal and $U^\transpose MU$ diagonal.
A simple computation shows that $U$ is symplectic.

If $d<n$, consider the orthogonal projection $P$
onto the span of $\iso u_1,\ldots,\iso u_d,u_1,\ldots,u_d$.
Then $P$ commutes with both $M$ and $\iso$.
So we can choose an orthonormal set of vectors
$\iso u_{d+1},\ldots,\iso u_n,u_{d+1},\ldots,u_n$
in the null space of $P$ and define $U$ as above.
\qed

\medskip
The same construction can be used to give a

\proofof(JoltEquiv)
The necessity of the condition $DX(x)DX(y)=0$ is obvious.
Assume now that this condition holds, for all $x,y\in\real^{2n}$,
and that $H$ is regular.

First, we show that $H$ is affine-integrable.
By regularity, there exist $x_0,w\in\real^n$ such that $DX(x_0)w=X(x_0)$.
Thus $X(x)=DX(x_0)w+\int_0^1 DX(x_0+sv)v\,ds$,
for any given $x\in\real^{2n}$, where $v=x-x_0$.
This shows that $DX(x)X(x)=0$ for all $x$,
which implies that $H$ is affine-integrable.

Let $(u_1,u_2,\ldots,u_d)$ be an orthonormal basis
for the subspace spanned by all vectors $\Hess(y)z$
with $y,z\in\real^{2n}$.
Then $\Hess(x)\iso u_j=0$ for all $x$ and all $j$.
Defining $U$ as in the proof of \clm(DXNormalForm),
we obtain \equ(DXNormalForm) simultaneously for all $x$.
(The matrix $A$ can depend on $x$ and need not be diagonal.)
So $(H\circ U)(q,p)$ is independent of $q$ by \clm(BZzero),
implying that $H$ is a shear Hamiltonian.
\qed

\smallskip
As a corollary we obtain

\claim Theorem(MoserQSymp) {\rm[\rMoser]}
Every affine-integrable Hamiltonian $H$ that is homogeneous of degree $3$ is a shear.

\proof
By \equ(CommDuHDvH) we have $\{\partial_u^k H,H\}=0$ for $k\le 2$.
The same holds for $k\ge 3$ since $H$ is of degree $3$.
It follows that $\{H(\bdot+u),H\}=0$ for all $u$.
Or equivalently, $X(x)^\transpose\iso X(y)=0$ for all $x$ and $y$.
{}From this we get $DX(x)DX(y)=0$ by differentiation,
and the assertion follows from \clm(JoltEquiv).
\qed

\medskip
The matrix $U$ described in \clm(DXNormalForm)
is both symplectic and orthogonal.
This means that $U^\transpose\iso U=\iso$ and $U^\transpose U=\id$.
As a result, we also have $\iso U=U\iso$.
In fact, any two of the three properties imply the third.
This is know as the $2$-out-of-$3$ property of the unitary group
$\rmU(n)=\rmO(n)\cap\Sp(2n,\real)\cap\GL(n,\complex)$.
The complex structure here is given by the matrix $\iso$,
and the equation $\iso U=U\iso$ simply says that $U$ is ``complex''.
Using the properties $U^\transpose U=\id$ and $\iso U=U\iso$,
any matrix $U\in\rmU(n)$ can be written as
$$
U=\twomat{S}{T}{-T}{S}\,,\qquad
S^\transpose S+T^\transpose T=\id\,,\quad
S^\transpose T=T^\transpose S\,.
\equation(UComplex)
$$
We will refer to such a $2n\times 2n$ matrix as being {\sl unitary}.
The $n\times n$ submatrices $S$ and $T$ will be referred to
as the real and imaginary parts of $U$, respectively.

\smallskip
Concerning the claim in \dem(Iwasawa), we note that
any symplectic matrix $M$ can be written as a product
$M=UAN$, where $U$ is unitary, $A$ positive diagonal,
and $N$ unipotent upper-triangular.
This is the standard Iwasawa decomposition [\rHelgason].
If $H$ is a Hamiltonian such that $(H\circ M)(q,p)$
is independent $q$, then $(H\circ U)(q,p)$
is independent of $q$ as well.

\smallskip
By \clm(DXNormalForm), the Hessian $\Hess(x)$
of an affine-integrable Hamiltonian $H$ is always of the form
$$
\Hess(x)=\twomat{S}{T}{-T}{S}\twomat{0}{0}{0}{A}
\twomat{S^\transpose}{-T^\transpose}{T^\transpose}{S^\transpose}
=\twomat{TAT^\transpose}{TAS^\transpose}
{SAT^\transpose}{SAS^\transpose}\,,
\equation(GeneralHess)
$$
for any given $x\in\real^{2n}$,
where the matrix $A$ can be chosen to be diagonal.
This representation is unique if $\Hess(x)$ has $n$
distinct nonzero eigenvalues, and if the diagonal elements of $A$
are required to be in some prescribed order.

If we do not require that the matrix $A$ be diagonal,
then we could replace $S$, $T$, and $A$ in \equ(GeneralHess)
by $SV$, $TV$ and $V^{-1}AV$, respectively, where $V$
can be any orthogonal $n\times n$ matrix.
This fact is used in the lemma below.

\demo Example(Trivial)
Let $M$ be an $m\times n$ matrix of rank $m\le n$,
and let $f:\real^m\to\real^n$ be real analytic.
Assuming $Mf=0$, the equation
$\dot q=f(Mq)$ defines a flow on $\real^n$
that is linear in time: $\Phi^t(q_0)=q_0+tf(Mq_0)$.
This is similar to the flow considered in [\rMcLachlan, Lemma 5].
It extends to a Hamiltonian flow on $\real^{2n}$,
as does every flow on $\real^n$.
The Hamiltonian is $H(q,p)=p^\transpose f(Mq)$.
Using \clm(JoltEquiv), is is easy to check that $H$ is a shear Hamiltonian.
In fact, $H$ can be trivialized explicitly:
If we set $T=M^\transpose(MM^\transpose)^{-1}M$ and $S=\id-T$,
then \equ(UComplex) defines a matrix $U\in\rmU(n)$,
and we get $(H\circ U)(q,p)=H(p,p)$.

\claim Definition(seminormal)
We say that $\Hess(x)$ is in semi-normal form
if $\ZZ(x)=0$ and $\BB(x)=0$.

\claim Lemma(Utheta)
Let $H$ be an affine-integrable Hamiltonian and $x\in\real^{2n}$.
If $\AA(x)$ is nonsingular then
$U^\transpose\Hess(x)U$ is in semi-normal form for the matrix
$$
U=\exp\twomat{0}{\theta}{-\theta}{0}\,,\qquad
\theta(x)=\tan^{-1}\bigl(\zeta(x)\bigr)\,,\qquad
\zeta(x)=\ZZ(x)\AA(x)^{-1}\,.
\equation(Utheta)
$$

\proof
Define $\zeta=\zeta(x)$ as above.
A comparison with \equ(GeneralHess) shows that $\zeta=TS^{-1}$.
The conditions in \equ(UComplex) on $S$ and $T$
imply that $\zeta$ is symmetric, and that 
$S^\transpose(\id+\zeta^2)S=\id$.
Since $SV=|S|$ for some orthogonal matrix $V$,
we can choose $S$ to be a positive definite symmetric matrix.
The choice is then unique: $S=(\id+\zeta^2)^{-1/2}$.
Setting $\theta=\tan^{-1}\bigl(\zeta)$ we obtain
$S=\cos(\theta)$ and $T=\zeta S=\sin(\theta)$,
which leads to the expression \equ(Utheta) for $U$.
\qed

This offers another way of checking whether $H$ is a shear Hamiltonian.
First, we note that a nondegenerate affine-integrable Hamiltonian $H$ is regular:
$X(x)$ belongs to the range of $DX(x)$
at every point $x$ where $DX(x)$ has rank $n$,
since $DX(x)X(x)=0$ by \clm(HXDXinvariant).

\claim Lemma(zeta)
Let $H$ be an affine-integrable Hamiltonian.
Assume that $\AA(x_0)$ is nonsingular at some point $x_0$.
Then $H$ is a shear if and only if $\zeta$ is constant near $x_0$.

\proof
First, assume that $\zeta$ is constant near $x_0$.
So near $x_0$, the matrix $U$ in \clm(Utheta) is independent of $x$,
and $H\circ U$ is in semi-normal form for a fixed
unitary matrix $U$.
By analyticity, this property extends to all $x\in\real^{2n}$.
Furthermore, $\Hess(x_0)$ has rank $n$, as \equ(GeneralHess) shows,
implying that $H$ is regular.
So $H$ is a shear Hamiltonian by \clm(BZzero).

Conversely, assume that $(H\circ U)(q,p)$ is independent of $q$
for some linear symplectic matrix $U$.
Then $U$ can in fact be chosen unitary,
as was shown the proof of \clm(JoltEquiv).
If $S$ and $T$ are the real
and imaginary parts of $U$, as defined by \equ(UComplex),
then we have $\ZZ(x)\AA^{-1}(x)=TS^{-1}$
at every point $x$ where $\AA(x)$ is nonsingular.
\qed

For completeness, let us mention that
there is an alternative representation of
$\Hess(x)$ via the shear map $(q,p)\mapsto(q,p+\zeta q)$.
Assuming that $H$ is affine-integrable and $\AA=\AA(x)$ nonsingular,
$$
\Hess(x)=\twomat{\id}{\zeta}{0}{\id}\twomat{0}{0}{0}{\AA}
\twomat{\id}{0}{\zeta}{\id}
=\twomat{\zeta\AA\zeta}{\zeta\AA}{\AA\zeta}{\AA}\,.
\equation(GeneralHessShear)
$$
This shear is not unitary.
But it is symplectic, since $\zeta$ is symmetric.
Furthermore, these shear maps form a group.
Notice also that $\BB(x)=\ZZ(x)\zeta(x)$.
So the condition $\BB(x)=0$ in \clm(seminormal) is redundant
if $H$ is affine-integrable and $\AA(x)$ nonsingular.

\section Invariant affine subspaces
%%%%%%%%%%%%%%%%%%%%%%%%%%%%%%%%%%%

In this section we give a proof of \clm(AffineSub)
and some related results.
It is always assumed that $H$ is affine-integrable and real analytic.

Besides the flow $\Phi$ for the Hamiltonian $H$,
consider also the flows $\Psi_j$ for the Hamiltonians
$\partial_jH$, where $\partial_jH$ denotes the $j$-th
partial derivative of $H$.
By standard ODE results, $\Psi_j^t(x)$ is well defined
for all times $t$ in some open neighborhood of zero in $\complex$
(which may depend on $x$).
By \equ(CommDuHDvH) the flows $\Psi_j$ commute with each other
and with $\Phi$.
So the flow $\Psi^w$ for $\partial_wH$ is given by
$$
\Psi^{wt}=\Psi_1^{tw_1}\circ\Psi_2^{tw_2}\circ\cdots\circ\Psi_{2n}^{tw_{2n}}\,.
\equation(Psiw)
$$
Again, $\Psi^w(x)$ is well defined for all $w$
in some open ball $B(x)\subset\complex^{2n}$ centered at the origin.
Furthermore, the ``group property'' 
$\Psi^u(\Psi^w(x))=\Psi^{u+w}(x)$ holds whenever
$w,w+u\in B(x)$ and $u\in B(\Psi^w(x))$.

\claim Lemma(DFZero)
Let $u\in\real^{2n}$.
If the derivative of $\partial_uH$ vanishes at some point $x$,
then it vanishes at $\Psi^w(x)$ for every $w\in B(x)$.

\proof
Define $(G)_t=G\circ\Psi^{tw}$
for any function $G$ on $\real^{2n}$.
Let now $G=\partial_uH$. Then
$$
\eqalign{
0&=\bigl(\partial_j\{G,\partial_wH\}\bigr)_t
=\bigl(\{\partial_jG,\partial_wH\}+\{G,\partial_j\partial_wH\}\bigr)_t\cr
&=\{(\partial_jG)_t,\partial_wH\}
+\{G,(\partial_j\partial_wH)_t\}\,.\cr}
\equation(DGZeroOne)
$$
Here, we have used that $\partial_wH$ and $G$ are invariant
under the flow $\Psi^w$, and that the maps $\Psi^{tw}$ are symplectic.
Thus, we have
$$
{d\over dt}(\partial_jG)_t
=\{(\partial_jG)_t,\partial_wH\}
=-\{G,(\partial_j\partial_wH)_t\}
=-\sum_{\sigma,\tau}(\partial_\sigma G)\iso_{\sigma,\tau}
\partial_\tau(\partial_j\partial_wH)_t\,.
\equation(DGZeroTwo)
$$
Due to the factors $\partial_\sigma G$ that all vanish at $x$,
the value $(\partial_jG)_t(x)$ is independent of $t$ and thus
$(\partial_jG)\bigl(\Psi^{tw}(x)\bigr)=(\partial_jG)_t(x)=(\partial_jG)(x)=0$,
for all $t$ in some open neighborhood of zero.
The assertion now follows from 
the above-mentioned group property of $\Psi$
and the analyticity of $G$.
\qed

\claim Corollary(NullInvar)
Let $u\in\real^n$.
If $\Hess u$ vanishes at some point $x$,
then $\Hess u$ vanishes at $\Psi^w(x)$ for every $w\in B(x)$.
In other words, the null space (and thus the range)
of $\Hess$ is invariant under $\Psi^w$.

Notice that the same holds for $DX=\iso\Hess$.

\claim Corollary(PsiAffine)
Let $x\in\real^{2n}$ and $R(x)=\range(DX(x))$.
Then $\Psi^w(x)$ belongs to the affine space
$x+R(x)$ for all $w\in B(x)$.
Furthermore, $w\mapsto\Psi^w(x)$ is locally (near zero)
invertible as a map from $\iso R(x)$ to $x+R(x)$.

\proof
Consider the curve $u(t)=\Psi^{tw}(x)-x$.
Clearly $u(0)$ belongs to $R(x)$.
The derivative $u'(t)=DX(\Psi^{tw}(x))w$
belongs to the range of $DX(\Psi^{tw}(x))$,
which agrees with $R(x)$ by \clm(NullInvar).
Thus, $u(t)$ belongs to $R(x)$ whenever $tw\in B(x)$.
Since $\Hess(x)$ is symmetric,
$DX(x)$ is invertible as a map from $\iso R(x)$ to $R(x)$.
Thus, by the implicit function theorem,
the same holds locally (near zero) for the map
$w\mapsto\Psi^w(x)-x$, whose derivative at $w=0$ is $DX(x)$.
\qed

\proofof(AffineSub)
For each $y$ in $x+R(x)$ there exists an open neighborhood
$B_y$ of $y$ in $x+R(x)$ that is included in the range of $f_y: w\mapsto\Psi^w(y)$.
This follows from \clm(PsiAffine).
The vector field $X$ is constant on each $B_y$
since each component $X_j$ is invariant under the flow $\Psi^w$.
Similarly for $H$.
Furthermore, the open sets $B_y$ cover the affine space $x+R(x)$,
and since this space is connected,
it follows that $X$ and $H$ are constant on $x+R(x)$.

Assume now that $DX(x)$ has rank $n$.
If $y=x+u+v$, with $u\in R(x)$ and $v\in R(x)^\ssperp$,
then
$$
X(y)=X(x)+DX(x+u)v+\OO\bigl(|v|^2\bigr)\,.
\equation(Xxuv)
$$
If $y\not=x$ is sufficiently close to $x$
then $|DX(x+u)v|$ is bounded from below
by a positive constant times $|v|$,
so we have $X(y)=X(x)$ if and only if $y-x=u\in R(x)$.
\qed

\section Proof of Theorem 1.8
%%%%%%%%%%%%%%%%%%%%%%%%%%%%%

We will write the given nilpotency condition on $DX$ in the form
$$
\Hess(x'')\iso\Hess(x)\iso\Hess(x')=0\,,\qquad\qquad
x'',x,x'\in\real^{2n}\,.
\equation(HuHvHwo)
$$
It is straightforward to check that this condition is necessary
for $H$ to be linear-symplectic\-ally conjugate
to a Hamiltonian of the form \equ(SimpleHam).

Assume now that $H$ is a nondegenerate (and thus regular)
real analytic affine-integra\-ble Hamiltonian that satisfies \equ(HuHvHwo).
Consider a point $x_0=(q_0,p_0)$ where $\Hess(x_0)$ has rank $n$.
By performing a unitary change of variables, if necessary,
we may assume that $\Hess(x_0)$ is in semi-normal form,
$$
\Hess(x_0)=\twomat{0}{0}{0}{\AA(x_0)}\,.
\equation(SemiNormal)
$$
In the case where $x'=x''=x_0$, the property \equ(HuHvHwo)
implies that $D_q^2H=0$, so that
$$
\Hess(x)=\twomat{0}{\ZZ(x)}{\ZZ(x)^\transpose}{\AA(x)}\,,
\qquad\qquad x\in\real^{2n}\,.
\equation(SpecialHess)
$$
Thus, our Hamiltonian $H$ has to be of the form
$$
H(q,p)=\KK(p)+q^\transpose\VV(p)\,,
\equation(mThreeHqp)
$$
with $\KK:\real^n\to\real$ and $\VV:\real^n\to\real^n$
real analytic.
Another consequence of \equ(HuHvHwo) is that the range of
$$
\iso\Hess(x)\iso\Hess(x')
=\twomat{\ZZ(x)^\transpose\ZZ(x')^\transpose}
{\ZZ(x)^\transpose\AA(x')-\AA(x)\ZZ(x')}{0}{\ZZ(x)\ZZ(x')}
\equation(HessxHessxp)
$$
is contained
in the null space of $\Hess(x'')$, for every $x''$.
If we take $x''=x_0$ then this implies that $\ZZ(x)\ZZ(x')$ vanishes.
But $\ZZ(q,p)=D\VV(p)$ and thus
$$
D\VV(p)D\VV(p')=0\,,\qquad\qquad p,p'\in\real^n\,.
\equation(DVpropOne)
$$
Finally, multiplying \equ(HessxHessxp) on the left by $\Hess(x'')$
and using that the result has to be the zero matrix,
we find that $\ZZ(x'')^\transpose\AA(x)\ZZ(x')$ vanishes.
In particular,
$$
D\VV(p'')^\transpose\Kess(p)D\VV(p')=0\,,
\qquad\qquad p'',p,p'\in\real^n\,,
\equation(DVpropTwo)
$$
where $\Kess=\nabla D\KK$ is the Hessian of $\KK$.

\medskip
Let $\RR$ be the linear span of all vectors $\VV(p)$, with $p\in\real^n$.
Or equivalently, $\RR$ is the linear span of all vectors
$D\VV(p)u$ with $p,u\in\real^n$.
Here we have used that $\VV(p_0)=0$.
Then \equ(DVpropOne) and \equ(DVpropTwo) imply that
$$
\VV(p+v)=\VV(p)\,,\qquad
\KK(p+v)=\KK(p)+D\KK(p)v\,,\qquad v\in\RR\,.
\equation(RRinvar)
$$
Let $d$ be the dimension of $\RR$.
If $d=0$ then there is nothing left to prove.
Consider now the case where $d>0$.
Since $\VV(p)=0$ for all $p\in\RR$, we also have $d<n$.
Let $\RR^\ssperp$ be the orthogonal complement of $\RR$ in $\real^n$.
Notice that $q^\transpose\VV(p)$ vanishes whenever $q\in\RR^\ssperp$.

Next we apply a unitary change of variables.
To simplify the description,
we rename current quantities by adding a subscript ``old''.
The change of variables is $(q_\ssold,p_\ssold)=(Sq,Sp)$,
with $S$ orthogonal, such that $\RR=S^{-1}\RR_\ssold$ is the span
of all vectors $P=(p_1,\ldots,p_d,0,\ldots,0)$.
Then $\RR^\ssperp$ is the span of all vectors
$\bar p=(0,\ldots,0,p_{d+1},\ldots,p_n)$.
And $\VV=S^{-1}\VV_\ssold S$ takes values in $\RR$.
{}From \equ(RRinvar) we see that $\VV(p)$ does not depend on $P$,
and that $\KK(p)$ is an affine function of $P$.
Setting $Q=(q_1,\ldots,q_d,0,\ldots,0)$,
the new Hamiltonian $H$ is of the form \equ(SimpleHam).
This concludes the proof of \clm(SimpleHam).

\section Commutators and quadratic functions
%%%%%%%%%%%%%%%%%%%%%%%%%%%%%%%%%%%%%%%%%%%%

The main goal in this section is to give a proof of \clm(SimplerHam)
on quartic Hamiltonians.
But some of the observations and computations
apply to other Hamiltonians as well.

The Hessian $\Hess(u)$ at any point $u\in\real^{2n}$
defines a quadratic function $x\mapsto\half x^\transpose\Hess(u)x$.
The flow generated by this function is linear in time, since $DX(u)^2=0$.
It is useful to know how these flows
for different vectors $u$ are related,
as the hypotheses in \clm(JoltEquiv) and \clm(SimpleHam) show.
What simplifies the situation for homogeneous quartic Hamiltonians
is that $x^\transpose\Hess(u)x=u^\transpose\Hess(x)u=\partial_u^2H(x)$,
and if $H$ is affine-integrable,
then $\partial_u^2H$ commutes with $H$ by \equ(CommDuHDvH).
This fact will be exploited below.

First we note that the Poisson bracket of two homogeneous quadratic functions
$$
F(x)=\thalf x^\transpose\FF x\,,\qquad
G(x)=\thalf x^\transpose\GG x\,,
\equation(FxGx)
$$
is again a homogeneous quadratic function,
$$
\{F,G\}(x)=\thalf x^\transpose\EE x\,,\qquad
\iso\EE=\thalf(\iso\FF)(\iso\GG)-\thalf(\iso\GG)(\iso\FF)\,.
\equation(KFG)
$$
Here, $\FF$, $\GG$, and $\EE$ are symmetric $2n\times 2n$ matrices.
The corresponding matrices
$\iso\FF$, $\iso\GG$, and $\iso\EE$ belong to $\sp(2n,\real)$.

We also need to compute some double commutators,
and not all functions involved are quadratic.
To simplify the expressions,
we use the operator notation $\llb H\rrb F=\{F,H\}$.
A straightforward computation shows that
$$
\eqalign{
\partial_u\partial_v\llb H\rrb
&=\llb\partial_vH\rrb\partial_u
+\llb\partial_uH\rrb\partial_v
+\llb H\rrb\partial_u\partial_v
+\llb\partial_u\partial_vH\rrb\cr
&=\partial_u\llb\partial_vH\rrb
+\llb\partial_uH\rrb\partial_v
+\llb H\rrb\partial_u\partial_v\,.\cr}
\equation(DDHComm)
$$
Let now $F$ be is a polynomial of degree $\le 2$.
Then $\partial_u\partial_vF$ is constant
and thus commutes with every function.
If in addition $F$ commutes with $H$,
then \equ(DDHComm) yields
$$
\eqalign{
-\llb\partial_u\partial_vH\rrb F
&=\llb\partial_vH\rrb\partial_u F
+\llb\partial_u H\rrb\partial_v F\,,\cr
0&=\partial_u\llb\partial_vH\rrb F
+\llb\partial_uH\rrb\partial_vF\,.\cr}
\equation(duvHCommFSimple)
$$
So far we have not used any properties of $H$
other than differentiability.

\proofof(SimplerHam)
Assume now that $H$ is homogeneous of degree $4$.
Let $F$ be a polynomial of degree $\le 2$ that commutes with $H$.
Then the second identity in \equ(duvHCommFSimple) implies that
$$
\llb\partial_{v'}H\rrb\partial_{u'}
\llb\partial_vH\rrb\partial_u F
=-\llb\partial_{v'}H\rrb
\llb\partial_{u'}H\rrb\partial_v\partial_u F=0\,.
\equation(PreHupvpHuvFo)
$$
Assume in addition that $H$ is affine-integrable.
Then $H$ satisfies \equ(CommDuHDvH).
By Jacobi's identity for the Poisson bracket,
if $F$ and $G$ commute with $H$, then so does $\llb G\rrb F$.
In particular, $\llb\partial_u\partial_vH\rrb F$
commutes with $H$.
Using \equ(duvHCommFSimple) and \equ(PreHupvpHuvFo), we find that
$$
\llb\partial_{u'}\partial_{v'}H\rrb
\llb\partial_u\partial_vH\rrb F
=\bigl(\llb\partial_{v'}H\rrb\partial_{u'}
+\llb\partial_{u'}H\rrb\partial_{v'}\bigr)
\bigl(\llb\partial_vH\rrb\partial_u
+\llb\partial_uH\rrb\partial_v\bigr)F=0\,.
\equation(HupvpHuvFo)
$$
Since $H$ is homogeneous of degree $4$,
the second derivatives of $H$ are homogeneous quadratic polynomials.
They commute with $H$ by \equ(CommDuHDvH).
So as a special case of \equ(HupvpHuvFo) we have
$$
\llb\partial_u^2H\rrb
\llb\partial_v^2H\rrb
\partial_u^2H=0\,.
\equation(HuHvHuo)
$$
Applying \equ(KFG) with $F=\partial_u^2H$ and $G=\partial_v^2H$ yields
$\iso\EE=2\iso\Hess(u)\iso\Hess(v)-2\iso\Hess(v)\iso\Hess(u)$.
According to \equ(HuHvHuo) we have
$\iso\EE\iso\Hess(u)-\iso\Hess(u)\iso\EE=0$, which implies that
$$
\Hess(u)\iso\Hess(v)\iso\Hess(u)=0\,,
\qquad\qquad u,v\in\real^{2n}\,.
\equation(HuHvHuo)
$$
Here we have used that $\Hess(u)\iso\Hess(u)=0$.
The identity \equ(HuHvHuo) also be written as
$$
D^3X(u,u,w)^\transpose\Hess(v)D^3X(u,u,z)=0\,,
\qquad\qquad u,v,w,z\in\real^{2n}\,.
\equation(PreHuHvHwo)
$$
By Lemma 5.1 below, the same holds if each of the
arguments $u$ is replaced by a different vector in $\real^{2n}$.
In particular, we have \equ(HuHvHwo)
and thus $DX(x'')DX(x)DX(x')=0$ for all $x'',x,x'\in\real^{2n}$.
The assertion now follows from \clm(SimpleHam).
\qed

\claim Lemma(ppuppv)
Let $X,Y,Z$ be vector spaces.
Let $\langle\ldots\rangle: X^3\to Y$ be a symmetric cubic form
and $\odot: Y^2\to Z$ be a symmetric quadratic form.
Assume that $\langle p,p,u\rangle\odot\langle p,p,v\rangle=0$
for all $p,u,v\in X$.
Then $\langle u_1,u_2,u_3\rangle\odot\langle v_1,v_2,v_3\rangle=0$
for any $u_i,v_j\in X$.

\proof
Under the given assumption we have
$\langle p,p,p\rangle\odot\langle p,p,p\rangle=0$.
Now ``differentiate'' this identity twice:
Replace $p$ by $p+u+v$, expand, and then collect all terms
that are bilinear in $(u,v)$. The result is
$$
12\langle p,u,v\rangle\odot\langle p,p,p\rangle
+18\langle p,p,u\rangle\odot\langle p,p,v\rangle=0\,.
\equation(CQOne)
$$
By assumption, the second term vanishes,
so $\langle p,u,v\rangle\odot\langle p,p,p\rangle=0$.
Differentiating this identity once we get
$$
\langle w,u,v\rangle\odot\langle p,p,p\rangle
+3\langle p,u,v\rangle\odot\langle p,p,w\rangle=0\,.
\equation(CQThree)
$$
Similarly, differentiating
$\langle p,p,u\rangle\odot\langle p,p,v\rangle=0$ once yields
$$
2\langle p,w,u\rangle\odot\langle p,p,v\rangle
+2\langle p,p,u\rangle\odot\langle p,w,v\rangle=0\,.
\equation(CQFour)
$$
Now \equ(CQThree) can be used to rewrite \equ(CQFour) as
$$
-{2\over 3}\langle v,w,u\rangle\odot\langle p,p,p\rangle
-{2\over 3}\langle p,p,p\rangle\odot\langle u,w,v\rangle=0\,.
\equation(CQFive)
$$
Or simplified, $\langle p,p,p\rangle\odot\langle u,v,w\rangle=0$.
The assertion now follows by polarization.
\qed

\section The example from Theorem 1.6
%%%%%%%%%%%%%%%%%%%%%%%%%%%%%%%%%%%%%

Consider the Hamiltonian \equ(SimpleHam),
with the variables $(Q,\bar q;P,\bar p)$ renamed to $x=(q,y;p,z)$,
$$
H(x)=K(z)+q^\transpose V(z)+p^\transpose W(z)\,,
\equation(rHKHd)
$$
where $q,p\in\real^d$ and $y,z\in\real^{n-d}$.
The corresponding vector field $X=(\dot q,\dot y;\dot p,\dot z)$ is given by
$\dot z=0$ and
$$
\dot q=W(z)\,,\quad
\dot p=-V(z)\,,\quad
\dot y_j=\partial_j K(z)
+q^\transpose\bigl[\partial_j V(z)\bigr]
+p^\transpose\bigl[\partial_j W(z)\bigr]\,.
\equation(doty)
$$
Since $\dot z=0$, the Hamiltonian $H$ is affine-integrable if and only if
$$
\eqalign{
\ddot y_j
&={\dot q}^\transpose\bigl[\partial_j V(z)\bigr]
+{\dot p}^\transpose\bigl[\partial_j W(z)\bigr]\cr
&=W(z)^\transpose\bigl[\partial_j V(z)\bigr]
-V(z)^\transpose\bigl[\partial_j W(z)\bigr]\cr}
\equation(ddoty)
$$
is equal to zero for all $j$.
This is precisely the condition \equ(WDVVDWo).
According to \clm(JoltEquiv), $H$ is a shear
if and only if $DX(x)DX(x')$ vanishes for all $x$ and $x'$.
Or equivalently, if and only if
$$
-X(x)^\transpose\iso X(x')
=W(z)^\transpose V(z')-V(z)^\transpose W(z')
\equation(XxXxp)
$$
vanishes for all $x$ and $x'$.

It should be noted that the case $d=1$ is trivial:
If $v(s)=V_1(x+s(x'-x))$ and $w(s)=W_1(x+s(x'-x))$
satisfy $wv'-vw'=0$,
then by the quotient rule of differentiation,
the functions $v$ and $w$ are constant multiples of each other.
So \equ(XxXxp) follows from \equ(ddoty).
In this case, $H(q,p)$ can be made either independent of $q$
via a change of variables $(q,p)\mapsto(q,p+cq)$,
or independent of $p$ via a change of variables $(q,p)\mapsto(q+cp,p)$.
Thus, $H$ is a shear Hamiltonian if $d=1$.

Notice also that, if the right hand side of \equ(ddoty)
is equal to zero, then it remains zero if $\partial_j$
is replaced by $\partial_j^2$.
Thus, the right hand side of \equ(XxXxp) is of the order $|z-z'|^3$.
This has motivated our choice of $V$ and $W$ below.

The Hamiltonian \equ(CounterExample) can be written as
$$
H(q,y;p,z)=q_1V_1(z)+q_2V_2(z)+p_1W_1(z)+p_2W_2(z)\,,
\equation(CounterExampleAgain)
$$
where $q,y,p,z\in\real^2$ and
$$
V_1(z)=z_1^3\,,\quad
W_1(z)=z_2^3\,,\quad
V_2(z)=\sqrt{3}z_1^2z_2\,,\quad
W_2(z)=-\sqrt{3}z_1z_2^2\,.
$$
Let us compute the right hand side of \equ(XxXxp),
with $z'$ replaced by $w$ in order to simplify notation.
If $w_2=z_2$ then
$$
\eqalign{
\bigl[W_1(z)V_1(w)+W_2(z)V_2(w)\bigr]
&-\bigl[V_1(z)W_1(w)+V_2(z)W_2(w)\bigr]\cr
&=z_2^3w_1^3-3z_1z_2^2w_1^2z_2-z_1^3z_2^3+3z_1^2z_2 w_1z_2^2\cr
&=\bigl(w_1^3-3z_1w_1^2-z_1^3+3z_1^2w_1\bigr)z_2^3\cr
&=(w_1-z_1)^3z_2^3\,.\cr}
\equation(rTi)
$$
This is clearly nonzero at some points, so $H$ cannot be a shear.
On the other hand,
$$
\bigl[W_1(z)\partial_jV_1(z)+W_2(z)\partial_jV_2(z)\bigr]
-\bigl[V_1(z)\partial_jW_1(z)+V_2(z)\partial_jW_2(z)\bigr]=0
\equation(rDj)
$$
holds for $j=1$, due to the factor $(w_1-z_1)^3$ in \equ(rTi).
By symmetry, we have an expression analogous
to \equ(rTi) if $w_1=z_1$, with the cubic factor being $(w_2-z_2)^3$.
So \equ(rDj) holds for $j=2$ as well.
Thus, $H$ is affine-integrable, as claimed in \clm(CounterExample).

\smallskip
A straightforward computation shows that
the Hamiltonian \equ(CounterExampleAgain) is nondegenerate:
The Hessian $\Hess(x)$ has rank $n=4$ whenever $z_1z_2\not=0$.
Another noteworthy fact is that the matrix $\zeta(x)$
defined in \equ(Utheta)
depends on $x$ only via the ratio $z_1/z_2$.
But the dependence is nontrivial, so by \clm(zeta),
this shows again that $H$ cannot be a shear Hamiltonian.

\section Elementary factorization
%%%%%%%%%%%%%%%%%%%%%%%%%%%%%%%%%

A classical theorem by Jung [\rJung] asserts that
the group (under composition) of polynomial automorphisms of the plane $\real^2$
is generated by affine automorphisms
and elementary shears $(q,p)\mapsto(q+s(p),p)$.
No general result of this type is known in dimensions higher than $2$.
\clm(CubicDecomp) covers the special case
of symplectic maps $F=\id+X$, with $X$ homogeneous of degree $3$.
Its proof is based on the following observation.

\claim Lemma(HamDecomp)
Let $H_0$ be a polynomial affine-integrable Hamiltonian
of the from \equ(SimpleHam). Write $H_0=H_1+H_2$, where
$H_1=Q^\transpose V(\bar p)$ and $H_2=K(\bar p)+P^\transpose W(\bar p)$.
Then $H_1$ and $H_2$ Poisson-commute with $H_3=\half\{H_1,H_2\}$.
Furthermore, $H_1$ and $H_4=H_2-H_3$ are shear Hamiltonians,
and the corresponding time-one maps satisfy
$$
\Phi_{H_0}^1=\Phi_{H_1}^1\circ\Phi_{H_4}^1\,.
\equation(HamDecomp)
$$

\proof
A straightforward calculation yields $H_3=\half V(\bar p)^\transpose W(\bar p)$.
So each of the Hamiltonians $H_j$ is of the form \equ(SimpleHam).
Here, and in what follows, $0\le j\le 4$.
In addition, $H_j$ satisfies the affine-integrability condition \equ(WDVVDWo).
Thus, the adjoint map $\llb H_j\rrb: G\mapsto\{G,H_j\}$
has the following nilpotency property:
If $f$ is any polynomial, then $\llb H_j\rrb^kf=0$
for sufficiently large $k$.
Furthermore, $H_j$ commutes with $H_3$,
since $H_3$ only depends on the variable $\bar p$,
while $H_j$ is independent of the variable $\bar q$.
Thus, by the Baker-Campbell-Hausdorff formula,
$$
e^{t\llb H_0\rrb}f=e^{t\llb H_1\rrb}e^{t\llb H_2\rrb}e^{-t^2\llb H_3\rrb}f
=e^{t\llb H_1\rrb}e^{t\llb H_2-tH_3\rrb}f\,,
\equation(BCH)
$$
for every polynomial $f$.
To be more precise, \equ(BCH) is an identity for formal power series.
But due to the above-mentioned nilpotency property,
only finitely many terms of the series are nonzero.
So \equ(BCH) holds as an identity between polynomials.
Using that $f\circ\Phi_{H_j}^t=e^{t\llb H_j\rrb}f$ for any polynomial $f$,
we obtain \equ(HamDecomp) from \equ(BCH).

Let $j\ge 1$. Then the vector field $X_j=\iso\nabla H_j$
satisfies $X_j(x)^\transpose\iso X_j(x')=0$ for all $x$ and all $x'$,
as can be seen from \equ(XxXxp).
This shows that $H_j$ is a shear Hamiltonian.
\qed

\demo Remark(ElemShearOneFour)
The time-one map for $H_4$ is an elementary shear,
$\Phi_{H_4}^1(q,p)=(q+\nabla h_4(p),p)$, where $h_4(p)=H_4(p,p)$.
The time-one map for $H_1$ is unitarily conjugate to an
elementary shear $S_1(q,p)=(q+\nabla h_1(p),p)$.
A straightforward computation, similar to the one in \dem(Trivial),
shows that $h_1(p)=H_1(p,p)$.

\proofof(CubicDecomp)
Let $F$ be a symplectic map on $\real^{2n}$ such that
$X=F-\id$ is a homogeneous polynomial of degree $m\ge 2$.
First, we prove that $X=\iso\nabla H$ for some affine-integrable Hamiltonian $H$,
using the same arguments as in [\rLomeli].
The symplecticity condition \equ(FSymplectic)
implies that $X$ satisfies the equation \equ(PhiIsoPhi) for $t=1$.
In this equation, the terms in square brackets have to vanish separately,
since they have different degrees of homogeneity.
The first of the resulting identities implies that
the derivative of $\iso X$ is a symmetric matrix.
Thus, by the Poincar\'e Lemma, $\iso X$ is the gradient of a function $-H$.
The second identity implies that $(DX)^2=0$.
Thus, $(DX)X=0$, since $X(x)=m^{-1}DX(x)x$ by homogeneity.
This shows that the flow for $X$ is linear in time.
In conclusion, $F$ is the time-one map of an affine-integrable
Hamiltonian $H$ that is homogeneous of degree $m+1$.

Consider now $m=3$, and assume that $X$ is nondegenerate.
By \clm(SimplerHam), there exists a linear symplectic
map $U$ on $\real^{2n}$, such that $H_0=H\circ U^{-1}$
is a Hamiltonian of the form \equ(SimpleHam).
In fact, $U$ is unitary, as seen in Section 4.
Using \clm(HamDecomp), we have
$$
F=\Phi_H^1=U^{-1}\circ\Phi_{H_0}^1\circ U
=U^{-1}\circ\Phi_{H_1}^1\circ\Phi_{H_4}^1\circ U
=\Phi_{H_1\circ U}^1\circ\Phi_{H_4\circ U}^1\,,
\equation(HamDecompTwo)
$$
where $H_1$ and $H_4$ (and thus $H_1\circ U$ and $H_4\circ U$)
are shear Hamiltonians.
\qed

%%%%%%%%%%%
\references
%%%%%%%%%%%

{\ninepoint

\item{[\rJung]}
H.W.E.~Jung,
{\it \"Uber ganze birationale Transformationen der Ebene},
J. Reine Angew. Math. {\bf 184} 161-174 (1942).

\item{[\rDragtFinn]}
A.J.~Dragt and J.M.~Finn,
{\it Lie series and invariant functions for analytic symplectic maps},
J. Math. Phys. {\bf 17} 2215-2227 (1976).

\item{[\rHelgason]} S.~Helgason,
{\it Differential Geometry, Lie Groups, and Symmetric Spaces},
Academic Press, New York (1978).

\item{[\rDragtAbell]}
A.J.~Dragt and D.T.~Abell,
{\it Jolt Factorization of Symplectic Maps},
Int. J. Mod. Phys. A (Proc. Suppl.) {\bf 2B}, 1019-1021 (1993).

\item{[\rMoser]}
J.~Moser,
{\it On quadratic symplectic mappings},
Math. Z. {\bf 216}, 417-430 (1994).

\item{[\rFengWang]}
K.~Feng and D.-L.~Wang,
{\it Variations on a theme by Euler},
J. Comput. Math. {\bf 16}, 97-106 (1998).

\item{[\rBE]}
M.~de Bondt, A.~van den Essen,
{\it Singular Hessians}, J.~Algebra {\bf 282}, 195-204 (2004).

\item{[\rLomeli]}
H.E.~Lomel\'\i,
{\it Symplectic homogeneous diffeomorphisms, Cremona maps and the jolt representation},
Nonlinearity {\bf 18}, 1065-1071 (2005).

\item{[\rForest]}
E.~Forest,
{\it Geometric integration for particle accelerators},
J. Phys. A: Math. Gen. {\bf 39}, 5321-5377 (2006)

\item{[\rBondt]}
M.~de Bondt,
{\it Quasi-translations and counterexamples to the Homogeneous Dependence Problem},
Proc. Amer. Math. Soc. {\bf 134}, 2849-2856 (2006).

\item{[\rMcLachlan]}
R.I.~McLachlan, H.Z.~Munthe-Kaas, G.R.W.~Quispel, A.~Zanna,
{\it Explicit Volume-Preserving Splitting Methods
for Linear and Quadratic Divergence-Free Vector Fields},
Found. Comput. Math. {\bf 8} 335-355 (2008).

\item{[\rKL]}
See reference [\rLomeli].
For completeness, a proof is included
in our proof of \clm(CubicDecomp) in Section 7.

}
\bye